\documentclass[11pt]{article}
\usepackage{multicol}
\usepackage{amssymb,amsmath,amsthm,bm}
\usepackage[colorlinks=true, urlcolor=blue,linkcolor=blue, citecolor=blue]{hyperref}
\usepackage{mathrsfs,verbatim}
\usepackage{caption,cleveref }
\usepackage{graphicx, float}
\usepackage{color, mathtools, tikz, xcolor}
\usepackage{enumerate,enumitem}
\usepackage{tikz}
\usepackage{lineno}
\usetikzlibrary{shapes.geometric}
\usepackage{pdfpages}
\usepackage{bbm}
\usepackage{algorithm}
\usepackage{algorithmic}
\usetikzlibrary{calc}

\usepackage[mathscr]{euscript}
\usepackage{appendix}
\usepackage{geometry}
\numberwithin{equation}{section}

\renewcommand{\l}{\left}
\renewcommand{\r}{\right}

\usepackage{svg}
\usepackage{multicol}
\usepackage{amssymb,amsmath,amsthm,bm}
\usepackage[colorlinks=true, urlcolor=blue,linkcolor=blue, citecolor=blue]{hyperref}

\newcounter{propcounter}

\numberwithin{equation}{section}

	\newtheorem{theorem}{Theorem}[section]
	
        \newtheorem{question}[theorem]{Question}
	
	\newtheorem{coro}[theorem]{Corollary}
	\newtheorem{conjecture}[theorem]{Conjecture}
	
	\newtheorem{proposition}[theorem]{Proposition}
	\newtheorem{lemma}[theorem]{Lemma}

        \theoremstyle{definition}
	\newtheorem{defn}[theorem]{Definition}

\setlist{nolistsep}

\title{Unbounded degree spanning hypertrees in Dirac hypergraphs}
\author{
Yaobin Chen\thanks{Shanghai Center for Mathematical Sciences,~Fudan University,~Shanghai,~200438,~China.~{\tt ybchen21@m.fudan.edu.cn}.
Supported by National Natural Science Foundation of China grant 123B2012.}
\and
Seonghyuk Im\thanks{Department of Mathematical Sciences, KAIST, South Korea and Extremal Combinatorics and Probability Group (ECOPRO), Institute for Basic Science (IBS), Daejeon, South Korea.~{\tt seonghyuk@kaist.ac.kr}  Supported by the National Research Foundation of Korea (NRF) grant funded by the Korea government(MSIT) No. RS-2023-00210430 and by the Institute for Basic Science (IBS-R029-C4).}
\and
Junchi Zhang\thanks{Shanghai Center for Mathematical Sciences,~Fudan University,~Shanghai,~200438,~China.~{\tt jczhang24@m.fudan.edu.cn}.
}
}
\date{\today}

\begin{document}

\maketitle

\begin{abstract}

In 2001, Koml{\'o}s, S{\'a}rk{\"o}zy, and Szemer{\'e}di proved that every sufficiently large $n$-vertex graph with minimum degree at least $\left(1/2+\gamma\right)n$ contains all spanning trees with maximum degree at most $cn/\log n$.
We extend this result to hypergraphs by considering \emph{loose hypertrees}, which are linear hypergraphs obtained by successively adding edges that share exactly one vertex with a previous edge.

For all $k > \ell \geq 2$, we determine asymptotically optimal $\ell$-degree conditions that ensure the existence of all rooted spanning loose hypertrees, without any degree condition, in terms of the $(\ell-1)$-degree threshold for the existence of a perfect matching in $(k-1)$-graphs.

As a corollary, we also asymptotically determine the $\ell$-degree threshold for the existence of bounded degree spanning loose hypertrees in $k$-graphs for $k/2 < \ell < k$, confirming a conjecture of Pehova and Petrova in this range.
In our proof, we avoid the use of Szemer{\'e}di's regularity lemma.

\end{abstract}

\section{Introduction}

The study of spanning trees in graphs has been an active area of research for many years. 
A classical result of Koml{\'o}s, S{\'a}rk{\"o}zy and Szemer{\'e}di~\cite{Tree} showed that every sufficiently large $n$-vertex graph with minimum degree at least $\left(\tfrac{1}{2}+\gamma\right)n$ contains all spanning trees of bounded degree, thus confirming a conjecture of Bollob{\'a}s. This result is also notable as one of the earliest applications of the blow-up lemma. 
In 2001, Koml{\'o}s, S{\'a}rk{\"o}zy, and Szemer{\'e}di~\cite{komlos2001spanning} relaxed the maximum degree condition by showing that for each $\gamma > 0$, every sufficiently large $n$-vertex graph with minimum degree at least $\left(\tfrac{1}{2}+\gamma\right)n$ contains every spanning tree with maximum degree smaller than $cn/\log n$. This bound is tight up to the constant $c$.
These results have also been generalised to various other settings, including random graphs, pseudorandom graphs, and directed graphs; see~\cite{Bottcher2019, csaba2010tight, han2022spanning, hyde2023spanning, Richard2022, krivelevich2010embedding, krivelevich2017bounded, Richardspanningtree2019, mycroft2020trees}.

In this paper, we consider the corresponding version of this result for hypergraphs and also discuss the maximum degree of spanning trees. A $k$-graph is a hypergraph with each hyperedge exactly $k$ vertices. A $k$-graph is said to be \emph{linear} if every pair of distinct edges shares at most one common vertex. A \emph{loose hypertree} (also known as a \emph{linear tree}) is a connected, linear $k$-graph in which any two vertices are connected by a unique path. 
Equivalently, a loose hypertree can be obtained by iteratively adding edges such that each new edge intersects the existing set of vertices in exactly one vertex. Note that a spanning loose hypertree can only exist when the number of vertices $n$ satisfies $n \equiv 1 \pmod{k-1}$. Throughout this paper, we always assume this divisibility condition holds.  

Extending a result of Koml{\'o}s, S{\'a}rk{\"o}zy, and Szemer{\'e}di to hypergraphs, Pehova and Petrova~\cite{pehova2024embedding} established the first result on the minimum vertex degree required for the existence of bounded degree loose hypertrees. The minimum $\ell$-degree $\delta_\ell(H)$ of a
$k$-graph $H$ is the minimum number of edges containing any given set of $\ell$ vertices.
\begin{theorem}[\cite{pehova2024embedding}]
Let $\varepsilon > 0$ and $\Delta \in \mathbb{N}$. There exists an $n_0$ such that every $3$-graph $G$ on $n \geq n_0$ vertices with $\delta_1(G) \geq \left(\frac{5}{9}+\varepsilon\right)\binom{n}{2}$ contains every spanning loose hypertree $T$ with $\Delta_1(T) \leq \Delta$.
\end{theorem}
The leading constant $\frac{5}{9}$ is optimal, since there exists a loose hypertree $T$ containing a perfect matching, and a $3$-graph $G$ on $n$ vertices with $\delta_1(G) \geq \left(\frac{5}{9}-o(1)\right)\binom{n}{2}$ that contains no perfect matching.
They conjectured that the existence of a perfect matching is the essential barrier for embedding bounded degree loose hypertrees. 
To formally state the conjecture of Pehova and Petrova, we first define the minimum degree threshold.

We denote the \emph{perfect matching} $(k,\ell)$-\emph{threshold} $\delta^{\mathrm{PM}}_{k,\ell}$ to be the infimum $\delta$ such that for any~$\gamma > 0$, there exists an integer~$n_0$ such that any $k$-graph $G$ on $n \geq n_0$ vertices with $k \mid n$ and $\delta_\ell(G) > \left(\delta + \gamma\right)\binom{n}{k - \ell}$ contains a perfect matching. 
This value has received significant attention over many years. For example, Pikhurko~\cite{pikhurko2008perfect} showed that $\delta^{\mathrm{PM}}_{k,\ell} = 1/2$ if $\ell \in [k/2, k-1]$ and Treglown and Zhao~\cite{andrew} proved the exact minimal degree result. H{\'a}n, Person and Schacht~\cite{han2009perfect} showed that $\delta^{\mathrm{PM}}_{3,1} = 5/9$. See~\cite{Yisurvey} for a survey on this topic.
We now state the conjecture of Pehova and Petrova.

\begin{conjecture}[\cite{pehova2024embedding}]\label{conj1}
    Let $\gamma>0, \Delta,k\in \mathbb{N}$ and $\ell\in [k-1]$. There exists an $n_0$ such that every $k$-graph $G$ on $n\geq n_0$ vertices with $\delta_\ell(G)\geq \l(\delta^{PM}_{k,\ell}+\gamma\r)\binom{n}{k-\ell}$ contains every spanning $k$-hypertree $T$ with $\Delta_1(T)\leq \Delta$.
\end{conjecture}

Pavez-Sign{\'e}, Sanhueza-Matamala, and Stein~\cite{stein2020tree} studied the minimum degree conditions for the existence of bounded degree spanning `tight' trees. Since every loose hypertree can be extended into a tight hypertree, their result confirms Conjecture~\ref{conj1} when $\ell = k - 1$. 
Recently, the first author and Lo~\cite{chen2025embedding} proved that the conjecture also holds when $\ell = k - 2$.
We remark that the first author and Lo also conjectured that their result can be extended to the `rooted' case; that is, the existence of an embedding of a spanning hypertree that maps the root vertex to a given vertex of~$G$.

In this paper, we aim to remove the bounded degree condition and determine the degree threshold for the containment of \emph{any} rooted spanning loose hypertree.
In the case of $2$-graphs, the problem is not meaningful, as the minimal vertex degree would need to exceed $n-1$ due to the example of a spanning star. 
By a similar construction, one can observe that the $1$-degree threshold for the existence of any spanning loose hypertree in a $k$-graph must be $1$.
Our main result determines the $\ell$-degree threshold for the existence of all rooted spanning loose hypertrees in a $k$-graph without any bounded degree condition for every $2 \leq \ell < k$.

\begin{theorem}\label{main theorem}
    Let $\gamma > 0$, $k \geq 3$, and $2 \leq \ell < k$. Then there exists an $n_0$ such that the following holds. 
    Let $G$ be a $k$-graph on $n \geq n_0$ vertices with $\delta_{\ell}(G) \geq \left( \delta^{\mathrm{PM}}_{k-1,\ell-1} + \gamma \right)\binom{n}{k - \ell}$ and let $v \in V(G)$. 
    Then for every loose hypertree $T$ on $n$ vertices with the root vertex $r \in V(T)$, there exists an embedding of $T$ into $G$ such that the image of $r$ is $v$.
\end{theorem}

We remark that this result is asymptotically optimal. Let $H$ be a $(k-1)$-graph on $n-1$ vertices with $\delta_{\ell-1}(H) \geq \left( \delta^{\mathrm{PM}}_{k-1,\ell-1} - o(1) \right)\binom{n-1}{k - \ell}$,
which exists by the definition of $\delta^{\mathrm{PM}}_{k-1,\ell-1}$. 
We construct a $k$-graph $G$ by adding a vertex $v$ to $V(H)$ and including all edges of the form $e \cup \{v\}$ for every $e \in E(H)$, together with all the possible $k$-edges contained in $V(H)$. Then $G$ does not contain a star centred at $v$, and $\delta_{\ell}(G) \geq \left( \delta^{\mathrm{PM}}_{k-1,\ell-1} - o(1) \right)\binom{n}{k - \ell}$.

As mentioned before, a result of Pikhurko~\cite{pikhurko2008perfect}  states that $\delta^{\mathrm{PM}}_{k-1,\ell-1} = \delta^{\mathrm{PM}}_{k,\ell} = 1/2$ when $k/2 < \ell < k$. Hence, Theorem~\ref{main theorem} implies that Conjecture~\ref{conj1} holds in these cases in a stronger form that the bounded degree condition can be removed. 

\begin{coro}
Let $\gamma > 0$, $k \geq 3$, and $k/2 < \ell < k$. Then there exists an $n_0$ such that for every $k$-graph $G$ on $n \geq n_0$ vertices with
$\delta_{\ell}(G) \geq \left( 1/2 + \gamma \right)\binom{n}{k - \ell}$ contains every spanning loose hypertree.

In particular, Conjecture~\ref{conj1} holds for all $k \geq 3$ and $k/2 < \ell < k$.
\end{coro}

We also note that when $k \geq 5$, the value $\ell = k - 2$ satisfies the inequality $k/2 < \ell < k$. Thus, our result also confirms the conjecture of the first author and Lo~\cite{chen2025embedding} on the rooted version for the case $k \geq 5$.

\subsection{Proof outline}
We aim to use the absorption method to prove Theorem~\ref{main theorem}. As a first step, we show that it is possible to find an almost spanning hypertree in $G$.
To achieve this, we use a decomposition of hypertrees which was originally introduced by Montgomery, Pokrovskiy, and Sudakov~\cite{Montgomery2020embedding} for $2$-graphs, and extended to loose hypertrees by Im, Kim, Lee, and Methuku~\cite{Im2024}. 
Their result states that for any loose hypertree $T$, there exists a sequence of hypertrees $T_0 \subseteq T_1 \subseteq \cdots \subseteq T_\ell = T$ satisfying the following properties:
\begin{enumerate}[label=(\arabic*)]
    \item $|T_0| = o(n)$;
    \item $\ell = O(1)$;
    \item $T_1$ is obtained by adding large stars to $T_0$;
    \item For each $i \geq 1$, $T_i$ is obtained either by adding $o(n)$ paths of length $3$ whose endpoints are in $T_{i-1}$, or by adding a matching in which each edge contains exactly one vertex from $T_{i-1}$.
\end{enumerate}
See Section~\ref{sec:structure_hypertree} for the precise definitions and the statement.
We embed each $T_i$ into $G$ sequentially. The base case $T_0$ can be embedded by using a simple greedy algorithm. When $T_i$ is obtained by adding $o(n)$ paths of length $3$, we apply the standard reservoir technique to embed these paths.
The non-trivial cases are embedding $T_1$ or embedding $T_i$ when $T_i$ is obtained by adding a matching. To address these, we use a rainbow matching result, which is applied by defining an appropriate graph system derived from link graphs.

The next step is to complete the embedding via absorption. We aim to use an absorption structure similar to that used in earlier work of Pehova and Petrova~\cite{pehova2024embedding}. However, their structure is only applicable when the hypertree contains suitable substructures of bounded maximum degree.
To address this limitation, we apply two distinct strategies for absorbing the remaining vertices depending on the number of leaves in the tree. A key observation is that any loose hypertree must contain either many leaves or many bare paths of length $6$.
If there are many leaves, we may choose a random subset and apply the earlier rainbow matching argument to absorb the remaining vertices. On the other hand, if there are many bare paths, then we can use these bare paths to embed vertex-disjoint absorbers, which allows us to apply the absorption structure of Pehova and Petrova~\cite{pehova2024embedding}.

\textbf{Paper Organisation:} In Section~\ref{sec:preliminaries}, we introduce the necessary notation and preliminary lemmas. In Section~\ref{sec:almost_embedding}, we prove that any almost spanning hypertree can be embedded. Section~\ref{sec:absorption} we prove the absorption lemma. In Section~\ref{sec:main_theorem_proof}, we combine all the ingredients to prove the main theorem. Finally, we conclude the paper with a discussion and some open problems in Section~\ref{sec:concluding_remark}.

\section{Preliminaries}\label{sec:preliminaries}
\subsection{Notations}
We use the notation $\ll$ to indicate a hierarchy between constants. 
More precisely, we say that a statement holds for $\alpha \ll \beta$ if there exists a non-decreasing function $f:(0, 1] \to (0, 1]$ such that the statement holds whenever $\alpha < f(\beta)$. 
We do not specify the function $f$ unless it is necessary.
We write $[n]$ to denote the set $\{1, 2, \dots, n\}$.

We define the \emph{normalised minimum $\ell$-degree} of $G$ as $\overline{\delta}_\ell(G) := \delta_\ell(G)/\binom{n}{k-\ell}$.
We write $G - H$ for the subgraph of $G$ obtained by removing the edges of $H$ from $G$, and $G \setminus H$ for the subgraph obtained by removing the vertices of $H$. 
For a vertex subset $U \subseteq V(G)$, the $k$-graph $G[U]$ denotes the subgraph of $G$ induced on $U$.
For an integer $t \in [k-1]$ and a set $A \in \binom{V(G)}{t}$ of $t$ vertices, we define the degree $d_U(A)$ as the number of edges in $G[U \cup A]$ that contain $A$. 
When $U = V(G)$, we simply write $d_G(A)$. 
Let $N_G(A)$ the set of $(k-t)$-tuples in $V(G) \setminus A$ that are contained in edges of $G$ together with $A$.

Note that it is straightforward to check that the minimum $1$-degree is at least the minimum $\ell$-degree.

\begin{proposition}
    Let $G$ be a $k$-graph and let $\ell \in [k-1]$. Then $\overline{\delta}_1(G) \geq \overline{\delta}_\ell(G)$.
\end{proposition}

For any set $S \subseteq V(G)$ with size less than $k$ and any $U \subseteq V(G) \setminus S$, the \emph{link graph} of $S$ in $U$, denoted by $L_U(S)$, is the $(k - |S|)$-graph on $U$ whose edges are all $(k - |S|)$-tuples $T$ such that $T \cup S \in E(G[U])$. When $U = V(G) \setminus S$, we simply write $L(S)$.

\subsection{Structure of hypertree}\label{sec:structure_hypertree}

In this section, we collect several definitions on hypertrees and some structural lemmas that will be used in the proof. 
A \emph{loose path} is a $k$-graph on $(k-1)\ell+1$ vertices $\{v_0,\dots,v_{(k-1)\ell}\}$ such that $\{v_{(k-1)i},\dots,v_{(k-1)(i+1)}\}$ is an edge for $i\in \{0, \dots, \ell-1\}$. 
Each of the two $(k-1)$-tuples $\{v_0,\dots, v_{k-2}\}$ and $\{v_{(k-1)(\ell-1)+1},\dots, v_{(k-1)\ell}\}$ in the first and the last edges of $P$ is called an \emph{end pairs} of $P$. 
A \emph{star} of size $D$ (or a \emph{$D$-star}) is a $k$-graph $H$ on $(k-1)D+1$ vertices such that $D$ edges share exactly one vertex.

A hypergraph $T$ is a \emph{loose hypertree} if and only if it is linear and there is a unique loose path between any pair of distinct vertices. 
Throughout the paper, we simply call a loose hypertree a \emph{hypertree}. 
A $u$--$v$ \emph{path} $P$ (or a path $P$ between $u$ and $v$) means a path $P$ with $u$ and $v$ in each end tuple of $P$. The vertices other than $u$ and $v$ in a $u$--$v$ path $P$ are called the \emph{interval} vertices of $P$.

In a hypertree $T$, a \emph{semi-bare path $P$} is a path such that edges in $E(T)\setminus E(P)$ are only incident to the vertices in its end pairs. A \emph{bare path $P$ between $u$ and $v$} is a $u$--$v$ path such that edges in $E(T)\setminus E(P)$ are only incident to $u$ and $v$.

A \emph{leaf vertex} of a hypertree $T$ is a vertex $v\in V(T)$ of degree one such that the edge $e$ containing $v$ contains at least $k-1$ vertices of degree one. 
The edges that contain a leaf are called \emph{leaf edges} of $T$. 
For a leaf edge $e$, the vertex in $e$ of degree at least two is called the \emph{parent} of $e$. 
We say a vertex is a parent if it is a parent of some leaf edge. 
If a matching $M$ consists of leaf edges of a hypertree $T$, we say that $M$ is a \emph{leaf matching} of $T$. 
For a leaf matching $M$, the set of leaves of $T$ in $M$ is called the \emph{leaf set} of $M$. We simply say that a vertex subset $X$ of $T$ is a \emph{matching leaf set} of $T$ if $X$ is the leaf set of a leaf matching $M$ in $T$.

The following lemma is almost the same as the lemma in Im, Kim, Lee and Methuku~\cite{Im2024} but we need a condition on the root vertex to discuss the embedding of rooted trees.
\begin{lemma}\label{lem:tree_split}
Let $D, n\geq 2$ be integers and let $0<\mu <1$. For any $k$-hypertree $T$ with at most $n$ edges together with two distinct root vertices $r_1, r_2 \in V(T)$, there exist integers $\ell\leq 10^5kD\mu^{-2}$ and $s\in [\ell-1]$ and a sequence of subgraphs $T_0\subseteq T_1\subseteq\dots \subseteq T_{\ell}=T$ such that the following holds:
    \stepcounter{propcounter}
\begin{enumerate}[label = {\bfseries \emph{\Alph{propcounter}\arabic{enumi}}}]
    \item $T_0$ has at most $\mu n$ edges, at most $k\mu n$ vertices, and $r_1, r_2 \in V(T_0)$; \label{cond:tree_split_1}
    \item $T_1$ is obtained by adding stars of size at least $D$ to $T_0$, that is, take pairwise vertex-disjoint stars of size at least $D$ and identify their centres with vertices in $T_0$; \label{cond:tree_split_2}
    \item For $i\notin \{0,s\}$, $T_{i+1}$ is obtained by adding a matching to $T_i$ such that $V(T_{i+1})\setminus V(T_i)$ is a matching leaf set of $T_{i+1}$; \label{cond:tree_split_3}
    \item $T_{s+1}$ is obtained by adding at most $\mu n$ vertex-disjoint bare paths of length $3$ to $T_s$ such that every bare path we add is a $u$--$v$ path $P$ where $u,v\in V(T_s)$, and $V(P)\setminus \{u,v\}$ is disjoint from $V(T_s)$. \label{cond:tree_split_4}
\end{enumerate}
\end{lemma}
For completeness, we state the proof of it in the appendix.
Recall that we use $H-G$ to denote the hypergraph obtained from $H$ by removing all edges in $G$.

\begin{lemma}[{\cite[Lemma~3.1]{Im2024}}]\label{count bare path}
Let $\ell,m\geq 2$ be integers. Let $T$ be a hypertree with at most $\ell$ leaf edges. Then there exist edge-disjoint semi-bare paths $P_1,P_2,\dots,P_s$ of length $m+1$ such that 

\begin{align*}
    e(T-P_1-\dots-P_s)\leq 6m\ell+\tfrac{2e(T)}{m+1}.
\end{align*}
\end{lemma}
For a tree $T$, it is well-known that $T$ has either many leaves or many bare paths of a fixed length (see, for example,~\cite{krivelevich2010embedding}). We generalise this fact to hypertrees as follows.
\begin{lemma}\label{lem:leaf_or_bare_path}
Let $1/n\ll \gamma, 1/k$. Then for a hypertree $T$ on $n$ vertices, there exist at least $\gamma n$ leaf edges or at least $2\gamma n$ semi-bare paths of length $6$.
\end{lemma}
\begin{proof}
    Suppose that there are at most $\gamma n$ leaf edges in $T$. Then by Lemma~\ref{count bare path} with $m=5$ and $\ell=2\gamma n$, we have
    \begin{align*}
        e(T-P_1-\dots -P_s)\leq 72\gamma n+\tfrac{2n}{6}.
    \end{align*}
    It follows that $e(P_1)+\dots +e(P_s)\geq 12\gamma n$. Hence there exist at least $2\gamma n$ semi-bare paths.
\end{proof}

The following lemma states that a tree on $n$ vertices contains a subtree of a desired small size.
\begin{lemma}[{\cite[Proposition 3.22]{Richardspanningtree2019}}]\label{lem:subtree_size}
    Let $n$ and $m$ be integers with $1 \leq m \leq n/3$. Let $T$ be a tree on $n$ vertices with root vertex $r \in V(T)$.
    Then there exist subtrees $T_1$ and $T_2$ such that $E(T_1)\cup E(T_2)=E(T)$, $|V(T_1)\cap V(T_2)|=1$, $r \in V(T_2)$, and $m \leq |V(T_1)| \leq 3m$.
\end{lemma}
We can extend the above lemma to hypertrees.

\begin{lemma}\label{lem:subtree_fixed_size}
Let $n,m,k\in \mathbb{N}$ and $2k \leq m\leq \tfrac{n}{3k}$. 
Let $T$ be a hypertree on $n$ vertices with root vertex $r \in V(T)$. Then there exist two edge-disjoint subtrees $T_1,T_2$ such that $V(T_1)\cup V(T_2)=V(T)$, $|V(T_1)\cap V(T_2)|=1$, $r \in V(T_2)$, and $\tfrac{m}{2k}\leq |V(T_1)| \leq 3m$.     
\end{lemma}
\begin{proof}
Let $L$ be the line graph of $T$. Let $T^*$ be the Breadth-First Search (BFS) tree of $L$ starting from any vertex.
Note that $|V(T^*)|=\tfrac{n-1}{k-1}$. 
Let $e_r$ be an arbitrary edge of $T$ containing $r$ and $m^* = \lceil \tfrac{m-1}{k-1}\rceil$.
By \Cref{lem:subtree_size}, there exist two vertex-disjoint subtrees $T_1^*$ and $T_2^*$ of $T^*$ such that $E(T_1^*)\cup E(T_2^*)=E(T^*)$, $|V(T_1^*)\cap V(T_2^*)|=1$, $e_r \in V(T_2^*)$, and $m^*\leq |V(T_1^*)| \leq 3m^*$. 
Let $\{e\} = V(T_1^*) \cap V(T_2^*)$, where $e$ is an edge of $T$ with $e=\{v_1,\dots,v_k\}$. 
Let $T_1'$ be the corresponding hypertree of $T_1^*$ in $T$. 
As $|V(T_1')| = (k-1)|V(T_1^*)|+1$, we have $m-k \leq |V(T_1')| \leq 3m$.
Then $T_1'-e$ is a forest with at most $k$ hypertrees. 
By the pigeonhole principle, there exists a component $T'$ of $T_1'-e$ with $|V(T_1')|/k \leq |V(T')| \leq |V(T_1')|$, so $m/2k \leq m/k-1 \leq |V(T')| \leq 3m$. 
Let $T_1=T'$ and let $T_2$ be the graph obtained by deleting all edges of $T_1$ from $T$. 
Then $E(T_1)\cup E(T_2)=E(T)$ and $|V(T_1)\cap V(T_2)|=1$. 
Furthermore, the edges of $T_2$ contain all the edges corresponding to the vertices in $V(T_2^*)$, so $r \in V(T_2)$ as required.
\end{proof}

\subsection{Concentration inequalities}

We will frequently use the following Chernoff bound.
\begin{lemma}[One-sided Chernoff Bound] Let $X_1,\dots ,X_n$ be mutually independent Bernoulli random variables and let $X=\sum_{i=1}^n X_i$. Then,
\begin{align*}
    &\mathbb{P}(X\geq (1+\varepsilon)\mathbb{E}[X])\leq \exp\left(-\frac{\varepsilon^2}{\varepsilon+2}\mathbb{E}[X]\right)\quad\text{ for every $\varepsilon>0$, and }\\
    &\mathbb{P}(X\leq (1-\varepsilon)\mathbb{E}[X])\leq \exp\left(-\frac{\varepsilon^2}{2}\mathbb{E}[X]\right)\quad \quad \text{ for every $\varepsilon\in (0,1)$ }.
\end{align*}
\end{lemma}

The following lemma by Kang, Kelly, K\"{u}hn, Osthus and Pfenninger~\cite{kang2023thresholds} will be used to prove degree concentration.
\begin{lemma}\label{lem:deg_concentration}
    Let $1/n \ll 1/s, \beta, \varepsilon<1$ and $k \geq 2$.
    Let $V$ be a set of size $n$ and let $p=p(n) \in [0, 1]$ satisfy $np \geq \varepsilon n^{\beta}$.
    Let $\mathcal{F} \subseteq \binom{V}{s}$ be a family of sets such that $|\mathcal{F}| \geq \varepsilon n^s (np)^{-1/2}$.
    Let $U \subseteq V$ be a $p$-random subset.
    Then with probability at least $1 - \exp(-n^{\beta/10})$, the number of elements of $\mathcal{F}$ contained in $U$ is $(1 \pm \varepsilon) p^s |\mathcal{F}|$. 
\end{lemma}

We now present a degree concentration lemma that we will frequently use in the proof.
\begin{lemma}\label{lem:deg_concentration2}
    Let $k, \ell$ be integers such that $k > \ell \geq 1$. Let $1/n \ll \gamma, \varepsilon, 1/k, 1/\ell <1$ and $\delta \in (0, 1)$.
    Let $G$ be a $k$-graph on $n$ vertices with $\overline{\delta}_\ell(G) \geq (\delta + \gamma)$.
    Let $p$ be a positive real number such that $\varepsilon \leq p \leq 1$ and $X$ be a $p$-random subset of $V(G)$.
    Then with probability at least $1-o(n^{-2})$, for every $\ell$-subset $S \subseteq V(G)$, we have 
    $$d_{X}(S) \geq (\delta + \gamma/2)\binom{np}{\ell}.$$
\end{lemma}
\begin{proof}
    We fix an $\ell$-subset $S \subseteq V(G)$ and let $\mathcal{F}_S \subseteq \binom{V(G)}{k-\ell}$ be the family of $(k-\ell)$-subsets $T$ such that $S \cup T \in E(G)$. 
    By the minimum degree condition, we have $|\mathcal{F}_S| \geq (\delta + \gamma)\binom{n}{k-\ell}$.
    Then by Lemma~\ref{lem:deg_concentration}, we have with probability at least $1-\exp(-n^{1/20})$, the number of elements of $\mathcal{F}_S$ contained in $X$ is at least $(1 - \gamma/4) p^{k-\ell} |\mathcal{F}_S| \geq (\delta + \gamma/2)\binom{np}{\ell}$.
    By taking a union bound over all $\ell$-subsets $S \subseteq V(G)$, we have with probability at least $1-o(n^{-2})$ that this inequality holds for every $\ell$-subset $S \subseteq V(G)$.
\end{proof}

\subsection{Rainbow perfect matching}
For a family of $k$-graphs $\mathbf{G}=(G_1, \dots, G_m)$ on the same vertex set $V$, we say that a graph on $H$ is called a \emph{rainbow subgraph} of $\mathbf{G}$ if there exists an injection $\phi: E(H) \to [m]$ such that $e \in E(G_{\phi(e)})$ for every edge $e \in E(H)$.
The study of rainbow subgraphs has received considerable attention in recent years. 
For example, Aharoni, Devos, de la Maza, Montejano and \v{S}\'{a}mal~\cite{Aharoni20rainbow} proved that if $e(G_i) \geq \lfloor \frac{26-2\sqrt{7}}{81}n^2\rfloor + O(n)$ for every $i \in [3]$, then $(G_1, G_2, G_3)$ contains a rainbow triangle.
Joos and Kim~\cite{Joos2020rainbow} proved that if $e(G_i) \geq n/2$ for all $i \in [n]$, then $(G_1, \dots, G_n)$ contains a rainbow Hamilton cycle.
See~\cite{Chakraborti2023bandwidth, Gupta2023general, Montgomery2022transversal} for more examples.
We will use the following theorem of Cheng, Han, Wang and Wang~\cite{rainbowspanning2023} regarding the existence of a rainbow perfect matching in a system of dense $k$-graphs.

\begin{lemma}\label{Jie Han}
    For every $\varepsilon>0$ and integer $d\in [k-1]$, there exists $n_0\in \mathbb{N}$ such that the following holds for all integers $n\geq n_0$ and $n\in k\mathbb{N}$. Every $n$-vertex $k$-graph system $\mathbf{G}=\{G_1,\dots, G_{m}\}$ with $m\leq n/k$ and $\overline{\delta}_d(G_i)\geq \delta^{PM}_{k,d}+\varepsilon$ for each $i\in [m]$ contains a rainbow matching with $m$ edges, i.e., disjoint edges $e_1,\dots, e_m$ with $e_j\in E(G_j)$ for $j\in [m]$.
\end{lemma}

\section{Almost embedding}\label{sec:almost_embedding}

In this section, we will prove the following lemma.
\begin{lemma}\label{lem:almost_spanning}
    Let $1/n \ll \mu, \gamma, 1/k$ and $2 \leq \ell < k$. Let $G$ be an $n$-vertex $k$-graph with $\overline{\delta}_{\ell}(G)\geq \left(\delta_{k-1, \ell-1}^{PM}+\gamma\right)$, and let $v_1, v_2 \in V(G)$ be two distinct vertices.
    Let $T$ be a hypertree on $(1-\mu)n$ vertices with two root vertices $r_1, r_2$.
    Then there exists an embedding $\phi: V(T) \to V(G)$ such that $\phi(r_1)=v_1$ and $\phi(r_2)=v_2$.
\end{lemma}

To prove this lemma, we use Lemma~\ref{lem:tree_split} and embed each $T_i$ greedily. 
The following lemma is a key tool for embedding large stars and matchings.
\begin{lemma}\label{lem:embed_stars}
    Let $k > \ell \geq 2$ and $1/n \ll \mu,\gamma, 1/k, 1/\ell$.
    Let $G$ be an $n$-vertex $k$-graph and let $X \subseteq V(G)$ be a set of size $m$ with $m\geq \mu n$ such that for every $\ell$-subset $S \subseteq V(G)$, we have $d_X(S) \geq (\delta_{k-1, \ell-1}^{PM}+\gamma)\binom{m}{\ell}$.
    Then the following holds.
    For every integer $t$, $t$ distinct vertices $v_1, \ldots, v_t \in V(G) \setminus X$ and integers $n_1, \ldots, n_t$ with $\sum_{i \in [t]} n_i \leq |X|/(k-1)$, there exist vertex-disjoint stars $S_1, \ldots, S_t$ centred at $v_1, \ldots, v_t$ respectively such that $|S_i| = n_i$ for all $i \in [t]$ and $V(S_i) \setminus \{v_i\} \subseteq X$. 
\end{lemma}
\begin{proof}
    Since we can enlarge the size of some star to adapt to the size of $X$, for simplicity, we may assume that $|X|$ is divisible by $k-1$ and $\sum_{i \in [t]} n_i = |X|/(k-1)$.
    For each $i \in [t]$, let $G_i$ be the link graph $L_X(\{v_i\})$.
    
    Then for every $S \subseteq X$ of size $\ell-1$, we have $d_X(S \cup \{v_i\}) \geq (\delta_{k-1, \ell-1}^{PM}+\gamma)\binom{|X|}{\ell}$.
    Thus, for each $i \in [t]$, we have $\delta_{\ell-1}(G_i) \geq (\delta_{k-1, \ell-1}^{PM}+\gamma) \binom{|X|}{\ell}$.
    We now define a system of graphs $\mathcal{H} = (H_1, \ldots, H_{|X|/(k-1)})$ on $X$ by taking $n_i$ copies of each $G_i$ for each $i \in [t]$.
    Then by \Cref{Jie Han}, there exists a rainbow perfect matching in $\mathcal{H}$.
    Thus, there exists a matching on $X$ such that it uses exactly $n_i$ edges from each $G_i$.
    Then by the definition of $G_i$, we can find vertex-disjoint stars $S_1, \ldots, S_t$ centred at $v_1, \ldots, v_t$ respectively such that $|S_i| = n_i$.
\end{proof}

We will also need to embed bare paths. The following lemma provides a key tool for this purpose.
\begin{lemma}\label{lem:bare_path}
    Let $2 \leq \ell < k$ and $\frac{1}{n} \ll \eta \ll \gamma, 1/k$. Let $G$ be an $n$-vertex $k$-graph with $\delta_{\ell}(G) \geq \gamma \binom{n}{k-\ell}$.
    Then for every two distinct vertices $u, v \in V(G)$, there exist at least $\eta n$ internally-disjoint paths of length three between $u$ and $v$.
\end{lemma}
\begin{proof}
    We choose a maximal collection $\mathcal{P}$ of internally-disjoint bare paths of length three between $u$ and $v$ and let $X$ be the set of vertices of the paths in~$\mathcal{P}$.
    If the number of paths in~$\mathcal{P}$ is less than $\eta n$, then we have $|X| \leq 3k\eta n$.
    Then for every $\ell$-subset $S$ of $V(G) \setminus X$, we have $d_{V(G) \setminus X}(S) \geq \gamma \binom{n}{k-\ell} - |X| \cdot n^{k-\ell-1} > \frac{\gamma}{2} \binom{n}{k-\ell}$.
    By taking $S$ be an arbitrary $\ell$-subset of $V(G) \setminus X \cup \{v\}$ that containing $u$, this shows that there exists an edge $e$ containing $u$ that avoids $X$ and $v$.
    Let $u'$ be an arbitrary vertex in $e \setminus \{u\}$.
    We now choose an arbitrary $(\ell-1)$-subset $W \subseteq V(G) \setminus (X \cup \{v\} \cup e)$.
    Then as $d_{V(G) \setminus (X \cup e)}(W \cup \{u'\}) \geq \frac{\gamma}{2} \binom{n}{k-\ell}$, we can choose an edge $f$ containing $W \cup \{u'\}$ that avoids $X$, $e$ and $v$.
    We choose a vertex $w \in W$ and let $W'$ be a $(\ell-2)$-subset in $V(G) \setminus (X \cup \{v\} \cup e \cup f)$.
    Then as $d_{V(G) \setminus (X \cup e \cup f)}(W' \cup \{w, v\}) \geq \frac{\gamma}{4} \binom{n}{k-\ell}$, we can choose an edge $g$ containing $W' \cup \{w, v\}$ that avoids $X$, $e$ and $f$.
    Then $e, f, g$ form a path of length three between $u$ and $v$, which contradicts the maximality of~$\mathcal{P}$.
\end{proof}

We are now ready to prove \Cref{lem:almost_spanning}.
\begin{proof}[Proof of \Cref{lem:almost_spanning}]
    Let $1/n \ll \alpha \ll \eta \ll \xi \ll \gamma, \mu, 1/k$.
    We first apply \Cref{lem:tree_split} with $D=3$ and $\eta$ plays the role of $\mu$ to obtain a sequence of subgraphs $T_0 \subseteq T_1 \subseteq \ldots \subseteq T_{\ell} = T$ that satisfies \ref{cond:tree_split_1}--\ref{cond:tree_split_4}.
    For each $i \in [\ell]$, we let $n_i = |V(T_i) \setminus V(T_{i-1})| + \mu n/3\ell$ and $p_i = n_i/n$.
    We also let $n_0=|V(T_0)| + \mu n/3, p_0=n_0/n$, $n_R = n-\sum_{i=0}^{\ell} n_i = \mu n/3$, and $p_R = n_R/n$.
    We partition $V(G) \setminus \{v_1, v_2\}$ into $V_0, \ldots, V_{\ell}$ randomly by selecting each vertex in $V(G)$ into $V_i$ with probability $p_i$ and into $R$ with probability $p_R$ independently.
    We claim that with high probability, this partition satisfies the following properties.
    \stepcounter{propcounter}
    \begin{enumerate}[label = {\bfseries \emph{\Alph{propcounter}\arabic{enumi}}}]
        \item $|V_i| = (1 \pm \alpha) n_i$ for every $i \in \{0\} \cup [\ell]$ and $|R| = (1 \pm \alpha) n_R$; \label{cond:random_partition_1}
        \item For every $\ell$-subset $S \subseteq V(G)$ and $i \in \{0\} \cup [\ell]$, we have $d_{V_i}(S) \geq \left(\delta_{k-1, \ell-1}^{PM}+\gamma/2\right) \binom{|V_i|}{\ell}$ and $d_R(S) \geq \left(\delta_{k-1, \ell-1}^{PM}+\gamma/2\right) \binom{n_R}{\ell}$; \label{cond:random_partition_2}
        \item For every pair of vertices $u \neq v \in V(G)$, there are at least $3k\eta n$ internally-disjoint paths of length three between $u$ and $v$ with all the internal vertices in $R$. \label{cond:random_partition_3}
    \end{enumerate}
    The first two properties follow from a standard application of the Chernoff bound and \Cref{lem:deg_concentration2}.
    For the last property, we fix two distinct vertices $u, v \in V(G)$.
    By \Cref{lem:bare_path}, there are at least $\xi n$ internally-disjoint paths of length three between $u$ and $v$.
    Let $\mathcal{P}$ be a collection of such paths.
    Then the events that all the internal vertices of a path $P \in \mathcal{P}$ are in $R$ are mutually independent and the probability of such an event is $p_R^{|V(P)|-2}$.
    By the Chernoff bound, we have with probability $1-o(n^{-2})$ that the number of paths in $\mathcal{P}$ with all the internal vertices in $R$ is at least $p_R^{|V(P)|-2} \xi n/2 \geq 3k\eta n$.
    By taking a union bound over all pairs of vertices $u, v \in V(G)$, we have with probability $1-o(1)$ that \ref{cond:random_partition_3} holds.

    We now fix a partition $V_0, \ldots, V_{\ell}, R$ of $V(G)$ that satisfies~\ref{cond:random_partition_1}--\ref{cond:random_partition_3}. 
    We now greedily embed $T_0, T_1, \ldots, T_{\ell}$ into $G$ by extending the previous embedding.

    \textbf{Embedding $T_0$:} We claim that we can embed $T_0$ into $G[V_0 \cup \{v_1, v_2\}]$ by a simple greedy algorithm.
    We first embed two roots $r_1$ and $r_2$ as $v_1$ and $v_2$, respectively.
    We order the edges of $E(T_0)$ by using BFS starting from $r_1$ to obtain a sequence $e_1, e_2, \ldots, e_m$. (Note that for components of $T_0$ not connected to $r_1$, we arbitrarily order them and start a BFS from an arbitrary vertex in each component.)
    We now embed the edges of $T_0$ one by one.
    Note that each $e_i$ contains at most one previously embedded vertex of $T_0$ except possibly when $e_i$ is the first edge containing $r_2$. In this case, $e_i$ contains at most two previously embedded vertices of $T_0$.

    Suppose that a subforest $F$ of $T_0$ has been embedded into $G[V_0 \cup \{v_1, v_2\}]$ and let $e$ be an edge of $E(T_0) \setminus E(F)$ such that $|e \cap V(F)| \leq 2$.
    Let $W$ be the set of vertices in the embedding of $F$.
    If $|e \cap V(F)| = 2$, then let $w_1, w_2 \in V_0$ be the images of the two vertices from $e \cap V(F)$. If $|e \cap V(F)| = 1$, then let $w_1 \in V_0$ be the image of the vertex from $e \cap V(F)$ and let $w_2$ be an arbitrary vertex in $V_0 \setminus W$.
    If $|V(e) \cap V(F)| = 0$, then let $w_1, w_2$ be two arbitrary unused vertices in $V_0 \setminus W$.
    
    As $|V_0| \geq \mu n/4$ and $\eta \ll \mu, 1/k$, the codegree of $w_1$ and $w_2$ in $V_0 \setminus W$ is 
    \begin{align*}
        d_{V_0 \setminus W}(w_1, w_2) &\geq \left(\delta_{k-1, \ell-1}^{PM}+\gamma/2\right) \binom{|V_0|} {k-2} - |W| \cdot |V_0|^{k-3} \\&\geq \delta_{k-1, \ell-1}^{PM}n_0^{k-2}/(k-2)! - k\eta n \cdot (1+\alpha)^{k-3}n_0^{k-3}> 0.
    \end{align*}

    Therefore, one can extend the embedding of~$F$ to $F \cup e$ by choosing an edge $e' \in E(G[V_0])$ containing $w_1$ and $w_2$ and whose other vertices are in $V_0 \setminus W$.
    By repeating this process, we can embed all edges in $T_0$ into $G[V_0 \cup \{v_1, v_2\}]$.

    \textbf{Embedding $T_i$ for $1 \leq i \leq s$:}
    We now assume that $T_{i-1}$ has been embedded into $\bigcup_{j=0}^{i-1} V_j$ and we want to extend the embedding to $T_i$.
    Let $u_1, \ldots, u_t$ be the vertices in $V(T_{i-1})$ that are incident to the edges in $E(T_i) \setminus E(T_{i-1})$.
    Let $v_1, \ldots, v_t$ be the images of $u_1, \ldots, u_t$ under the embedding in $\bigcup_{j=0}^{i-1} V_j$. 
    Let $m_j$ be the number of edges in $E(T_i) \setminus E(T_{i-1})$ incident to $u_j$.
    By the choice of $n_i$, we have $\sum_{j \in [t]} m_j = |V(T_i) \setminus V(T_{i-1})|/(k-1) \leq |V_i|/(k-1)$.
    By \ref{cond:random_partition_2}, we can apply \Cref{lem:embed_stars} with $t$, $v_1, \ldots, v_t$, and $m_1, \ldots, m_t$ to obtain vertex-disjoint stars $S_1, \ldots, S_t$ centred at $v_1, \ldots, v_t$ respectively such that $|S_j| = m_j$ for all $j \in [t]$ and $V(S_j) \setminus \{v_j\} \subseteq V_i$.
    Then we can extend the embedding of $T_{i-1}$ to $T_i$ by adding the stars $S_1, \ldots, S_t$.

    \textbf{Embedding $T_{s+1}$:} We now assume that $T_{s}$ has been embedded. 
    Let $v_1, \ldots, v_t$ and $u_1, \ldots, u_t$ be the vertices of $T_{s}$ such that adding the length-three bare paths between $u_i$ and $v_i$ for all $i \in [t]$ produces $T_{s+1}$.
    Let $v'_1, \ldots, v'_t$ and $u'_1, \ldots, u'_t$ be the embedded images of $v_1, \ldots, v_t$ and $u_1, \ldots, u_t$, respectively.
    By \ref{cond:random_partition_3}, for each pair $(u'_i, v'_i)$, we can find at least $3k\eta n$ internally-disjoint bare paths of length three.
    As $t \leq \eta n$, we can greedily choose these paths for each pair so that they are internally-disjoint from all previously chosen paths. This allows us to find $t$ such paths.
    Then we can extend the embedding of $T_{s}$ to $T_{s+1}$ by adding these paths.

    \textbf{Embedding $T_i$ for $s+1 < i \leq \ell$:} We follow the same strategy as for embedding $T_i$ for $1 \leq i \leq s$.
    By \ref{cond:random_partition_2} and \Cref{lem:embed_stars}, we can find a matching leaf set whose vertices are contained within $V_i$ so that we can extend the embedding of $T_{i-1}$ to $T_i$.
\end{proof}

\section{Absorption}\label{sec:absorption}
In this section, we will prove the absorption lemma as follows.
\begin{lemma}\label{absorb}
     Let $1/n\ll \varepsilon\ll \mu \ll \gamma, 1/k$ and $2 \leq \ell < k$. Let $G$ be a $k$-graph with $\overline{\delta}_{\ell}(G)\geq \left(\delta_{k-1, \ell-1}^{PM}+\gamma\right)$, and let $v\in V(G)$ be a vertex. Let $T$ be a hypertree on $\mu n$ vertices rooted at vertex $r$.

    Then there exists a vertex set $A\subseteq V(G)$ with $|A|=(\mu-\varepsilon)n$ and $v\in A$ such that for any $B\subseteq V(G)\setminus A$ with $|B|=\varepsilon n$, there exists an embedding $\psi$ of $T$ in $A\cup B$ with $\psi(r)=v$. 
\end{lemma}

We recall that we separate the proof of \Cref{absorb} into two cases based on the number of bare paths in $T$.
Although the absorber structure ``switcher'' as in \cite{chen2025embedding, pehova2024embedding, stein2020tree} is only used in the case of bounded degree, we will show that it still works when the number of bare paths is large.

\subsection{Many bare paths}
In this section, our aim is to prove the following absorption lemma. 
\begin{lemma}\label{absorb:many bare paths}
    Let $1/n\ll \varepsilon\ll \mu \ll \gamma$. Let $G$ be a $k$-graph with $\overline{\delta}_{1}(G)\geq \left(\frac{1}{2}+\gamma\right)$, and let $v\in V(G)$ be a vertex. Let $T$ be a hypertree on $\mu n$ vertices rooted at vertex $r$ with at least $2\mu^2 n$ edge-disjoint semi-bare paths of length~$6$.

    Then there exists a vertex set $A\subseteq V(G)$ with $|A|=(\mu-\varepsilon)n$ and $v\in A$ such that for any $B\subseteq V(G)\setminus A$ with $|B|=\varepsilon n$, there exists an embedding $\psi$ of $T$ in $A\cup B$ with $\psi(r)=v$. 
\end{lemma}

Note that $\delta^{\mathrm{PM}}_{k, \ell} \geq 1/2$ for every $k \geq 2$ and $\ell \in [k-1]$ so this implies \Cref{absorb} when $T$ has many bare paths.
To prove the lemma, we need to define an absorption structure.
For a star $S$ of size $2$, the \emph{leaf set} of $S$ is the set $N_S(v)$, which is a $(k-1)$-uniform matching of size $2$.
\begin{defn}\label{absorber construction}
    Let $G$ be a $k$-graph. For a $k$-ordered vertex-tuple $(w_1,\dots ,w_k)$, an \emph{absorbing tuple} for $(w_1, \dots, w_k)$ consists of $k-1$ vertex-disjoint $2$-stars $S_{v_2},\dots ,S_{v_k}$ centred at $v_2, \ldots, v_k$ respectively such that $w_1v_2\dots v_k\in E(G)$ and for each $2 \leq i \leq k$, there is a $2$-star $S_{w_i}$ with centre $w_i$ and $N_{S_{w_i}}(w_i)=N_{S_{v_i}}(v_i)$. 
    In other words, we can replace the centre of each $S_{v_i}$ with $w_i$ for each $2 \leq i \leq k$.
\end{defn}
Let $T$ be a $k$-hypertree with a root $r \in V(T)$ and $\psi:V(T)\rightarrow V(G)$ be an embedding of $T$ in $G$. 
We say that a $2$-star~$S_v$ centred at $v$ in $G$ is \emph{immersed} by $\psi$ if there exists a semi-bare path $P=e_1\dots e_6\subseteq E(T)$ with length~$6$ which does not contain $r$ as an internal vertex such that $v=\psi(e_3\cap e_4)$ and $e_3,e_4$ are mapped to distinct edges of $S_v$. 
We say that an absorbing tuple $(S_{v_2},\dots ,S_{v_k})$ is immersed by~$\psi$ if $S_{v_2},\dots, S_{v_k}$ are immersed by $\psi$.
For every distinct vertices $w_1,\dots,w_k\in V(G)$, let $A(w_1,\dots,w_k)$ be the family of absorbing tuples for $(w_1,\dots, w_k)$.

The next lemma shows one can extend an embedding of a hypertree using immersed absorbing tuples.

\begin{lemma}[Extension lemma]\label{Absorption lemma}
    Let $1/n\ll \mu \ll 1/k$. Let $G$ be a $k$-graph on $n$ vertices and $v \in V(G)$ be a fixed vertex. Let $T'$ be a $k$-hypertree on $\mu n-p(k-1)$ vertices with a root $r \in V(T')$. 
    Suppose that there exists an embedding $\psi$ of~$T'$ in $G$ with $\psi(T')=A$ and $\psi(r)=v$ and a family $\mathcal{A}$ of pairwise vertex-disjoint absorbing tuples such that 
    \stepcounter{propcounter}
    \begin{enumerate}[label = {\bfseries \emph{\Alph{propcounter}\arabic{enumi}}}]
        \item\label{absorb1} every absorbing tuple of $\mathcal{A}$ is immersed by $\psi$;
        \item\label{absorb2} for every distinct $w_1,\dots,w_k\in V(G)$, $\mathcal{A}$ contains at least $p$ absorbing tuples for $(w_1,\dots,w_k)$.
    \end{enumerate}
    Then for any $k$-hypertree $T$ on $\mu n$ vertices containing $T'$ and any $B\subseteq V(G) \setminus A$ with $|B|=p(k-1)$, there is an embedding $\varphi$ of $T$ in $A\cup B$ such that $\varphi(r)=\psi(r)$.
\end{lemma}

\begin{proof}
    We proceed by induction on $p$. The lemma holds when $p=0$, since $T=T'$. Thus, we assume that~$p\geq 1$.

    Consider any hypertree $T$ on $\mu n$ vertices containing $T'$. Let $e=\{x_1, \dots, x_k\} \in E(T)\setminus E(T')$ be an edge with $x_1\in V(T')$ and $x_2,\dots, x_k\notin V(T')$. 
    Let $T''=T'\cup \{x_1, \dots, x_k\}$. Let $u_1=\psi(x_1)$ and $u_2,\dots, u_k$ be distinct vertices in $B$. 
    Fix an absorbing tuple $\mathcal{S}=(S_{v_2},\dots, S_{v_k})$ for $(u_1,\dots,u_k)$ in~$\mathcal{A}$. We define $\psi':V(T'')\rightarrow V(G)$ as follows:

    \[
    \psi'(x) =
    \begin{cases}
    v_i &\quad \text{ if }x=x_i\text{ for some }2 \leq i \leq k, \\
    u_i &\quad \text{ if }x=\psi^{-1}(v_i) \text{ for some } 2 \leq i \leq k,\\
    \psi(x)&\quad \text{ otherwise}.
    \end{cases}
    \]
    We deduce that $\psi'$ is an embedding of $T''$ into~$G$. In addition, as $\psi(r)$ is not one of the centres of $S_{v_2},\dots,S_{v_k}$, we have $\psi'(r)=\psi(r)$.
    Let $\mathcal{A}'=\mathcal{A}\setminus \mathcal{S}$. 
    Then we observe that \ref{absorb1} and \ref{absorb2} are satisfied with $(T'',p-1,\psi',\mathcal{A}')$ playing the role of $(T',p,\psi,\mathcal{A})$.
    Since $T$ contains $T''$, our induction hypothesis implies that $G$ contains an embedding of $T$ into $A \cup B$ while preserving the image of the root $r$.
\end{proof}

The following lemma allows us to find many vertex-disjoint absorbing tuples when $\overline{\delta}_{1}(G)\geq 1/2+\gamma$.
\begin{lemma}[{\cite[Lemma 6.4]{chen2025embedding}}]\label{lemma:intersect}
    Let ${1}/{n}\ll \zeta \ll \alpha \ll \beta \ll \gamma \ll 1/k$. Let $G$ be a $k$-graph on $n$ vertices with $\overline{\delta}_{1}(G)\geq 1/2+\gamma$. Then there exists a set $\mathcal{A}$ of vertex-disjoint $2$-star tuples $\{S_{v_2^i},\dots ,S_{v_k^i}\}$ with $i\leq \beta n$ such that for distinct $w_1,\dots,w_k \in V(G)$, we have $|A(w_1,\dots ,w_k)\cap \mathcal{A}|\geq \alpha n$.
\end{lemma}

Next, we show that the absorbers can be immersed by a small tree if there are many bare paths.
\begin{lemma}[Immersing Lemma]\label{robust covering lemma}
    Let $1/n \ll \beta \ll \mu \ll \gamma\ll 1/k$. Let $G$ be a $k$-graph on $n$ vertices with $\overline{\delta}_{1}(G)\geq 1/2+\gamma$ and $T$ be a $k$-hypertree rooted at $r$ on $\mu n$ vertices with at least $\mu^2 n$ edge-disjoint semi-bare paths of length $6$. 
    Let $\mathcal{A}$ be a set of vertex-disjoint $2$-stars in $G$ with $|\mathcal{A}|\leq k \beta n$ and $v\notin V(\mathcal{A})$. 
    Then there exists an embedding $\psi: V(T)\rightarrow V(G)$ such that $\psi(r)=v$ and every star in $\mathcal{A}$ is immersed by $\psi$. 
\end{lemma}

\begin{proof}
    Let $\mathcal{A}=\{S_1, \dots, S_m\}$, where $m\leq k\beta n$. Let $T_0=\{r\}$ and $\psi_0:V(T_0)\rightarrow \{v\}$. 
    Let $P_1,\dots, P_\ell$ be the pairwise edge-disjoint semi-bare paths of length $6$ in $T$ which do not contain $r$ as an internal vertex. Then $\ell\geq \mu^2 n-1$. 
    Let $L$ be the line graph of $T$; i.e., $V(L)=E(T)$ and two vertices in $L$ are adjacent if they share a common vertex in $T$.
    We do a Depth-First Search (DFS) on $L$ starting at any edge $e$ containing $r$ and once it reaches an edge of $P_i$ for some $i \in [\ell]$, it will continue to search all the edges of $P_i$ before searching other edges.
    Let $e_1, \ldots, e_s$ be the resulting order of the DFS.
    
    For each $i \in [\ell]$, we choose the smallest index $t_i$ such that $e_1, \ldots, e_{t_i}$ contain exactly $i$ semi-bare paths among $P_1, \ldots, P_\ell$. 
    Let $T_i$ be the subtree of $T$ induced by $e_1, \ldots, e_{t_i}$.
    Without loss of generality, suppose that $T_i$ contains $P_1,\dots, P_{i}$. 
    Note that $T_i$ is edge-disjoint from all $P_j$ for $j<i$.
    We now claim that for each $i\in [\ell]$, there is an embedding $\psi_i:V(T_i)\rightarrow V(G)$ such that $S_j$ is immersed by $\psi_{i}$ for all $j \in [i]$ and $\psi_i(r)=v$.

    We proceed by induction on $i$. The base case $i=0$ is trivial.
    Suppose that the claim holds for $i-1$.
    Let $P_i=e^i_1\dots e^i_6$ and $x$ be the unique vertex in $e^i_3\cap e^i_4$.
    Since the semi-bare path $P_i$ is contained in $T_i$ and edge-disjoint from~$T_{i-1}$, the path from~$x$ to~$T_{i-1}$ has length at least $3$. Let $P'$ be the path from~$T_{i-1}$ to~$x$ in~$T_i$. Without loss of generality, suppose that $\{e_1^i,e_2^i,e_3^i\}= E(P_i \cap P')$. 
    Let $T'_{i}$ be the union of $T_{i-1}$, $P'$ and $e^i_4$. 
    As $P_i$ is a semi-bare path in $T$, all the edges incident to $x$ in $T$ are $e_3^i$ and $e_4^i$, and they are contained in $T_i'$. Moreover, this implies $T_i'\subseteq T_i$.

    Let $w$ be the centre vertex of $S_{i}$ and let the two edges of $S_i$ be $e_1 = \{w, u_2^1, \dots, u_k^1\}$ and $e_2 = \{w, u_2^2, \dots, u_k^2\}$. 
    We first greedily embed the edges of $P' \setminus P_i$. 
    As we embed at most $\mu n$ vertices, all the vertices $v \in V(G)$ incident to an edge that does not meet the previously embedded image except $v$. Thus, this greedy embedding is possible. 
    Then let the embedding function be $\psi_{i-1}'$. 
    Let $z$ be the vertex in $P_i\cap (P'-P_i)$. 
    Then $z\in e_1^i\setminus e_2^i$ by the definition of a semi-bare path.
    We now claim that there exists a path $e'_1e'_2$ from~$\psi_{i-1}'(z)$ to~$u_k^1$ in $G$ such that it is disjoint from the image of $\psi_{i-1}'$ except at $\psi_{i-1}'(z)$.
    Let $X$ be the image of $\psi_{i-1}'$.
    As $\overline{\delta}_{1}(G)\geq 1/2+\gamma$, we have that $d_{G \setminus X}(u_k^1) \geq \left(\frac{1}{2}+\frac{\gamma}{2}\right)\binom{n}{k-1}$ and $d_{G \setminus X}(\psi_{i-1}'(z)) \geq \left(\frac{1}{2}+\frac{\gamma}{2}\right)\binom{n}{k-1}$.
    Thus, the intersection of two link graphs $L_{G \setminus X}(\psi_{i-1}'(z))$ and $L_{G \setminus X}(u_k^1)$ has at least $\gamma \binom{n}{k-1}$ edges.
    Thus, it contains a star of size $2$ with edges $f_1$ and $f_2$. By taking $e'_1 = f_1 \cup \{\psi_{i-1}'(z)\}$ and $e'_2 = f_2 \cup \{u_k^1\}$, we have a desired path $e'_1e'_2$ in $G$.
    We now extend $\psi_{i-1}'$ by mapping $e_1^i$ to $e'_1$ and $e_2^i$ to $e'_2$.
    Then we can embed $e_3^i$ and $e_4^i$ onto the edges of $S_i$ such that $x$ is mapped to $w$. 
    For other edges in $T_{i}$, we can greedily embed into $G$ by the same argument as before. 
    Let $\psi_i$ be the embedding function of $T_i$.
    Then $S_1, \ldots, S_i$ are immersed by $\psi_i$ since $S_i$ is a $2$-star and the edges $e_3^i$ and $e_4^i$ are mapped to distinct edges of $S_i$. 
    Additionally, $\psi_i(r)=v$ as $\psi_i$ is an extension of $\psi_{i-1}$.

    By induction, we have an embedding $\psi_m:V(T_m)\rightarrow V(G)$. 
    Then by greedily extending $\psi_m$ to $T$, we can embed all edges in $T$ into $G$, which is the desired embedding of $T$ in $G$.
\end{proof}

\begin{proof}[Proof of Lemma~\ref{absorb:many bare paths}]
    Let $1/n\ll \varepsilon \ll \beta \ll \mu$. 
    By Lemma~\ref{lemma:intersect}, there exists a set $\mathcal{A}$ of vertex-disjoint absorbing tuples $\{S_{v_2^i},\dots ,S_{v_k^i}\}$ with $i\leq \beta n$ such that $v$ is not contained in any of $2$-stars and for distinct $w_1,\dots,w_k \in V(G)$, we have 
    \[
    |A(w_1,\dots ,w_k)\cap \mathcal{A}|\geq \frac{\varepsilon}{k-1} n.
    \]

    Let $T'$ be a hypertree obtained by iteratively removing $\varepsilon n/(k-1)$ leaf edges from $T$. Note that every hypertree has at least two leaf edges, thus we can remove a leaf edge without removing the root. 
    Then as $\varepsilon \ll \mu^2$, there still exist at least $\mu^2 n$ bare paths of length $6$ in $T'$. 
    By Lemma~\ref{robust covering lemma}, there exists an embedding $\psi':V(T')\rightarrow V(G)$ such that $\psi'(r)=v$ and every absorbing tuple in $\mathcal{A}$ is immersed by $\psi'$. 
    Let $A=\psi'(V(T'))$.
    Then by Lemma~\ref{Absorption lemma} with $p=\frac{\varepsilon}{k-1}$, for every $B \subseteq V(G) \setminus A$ with $|B|=\varepsilon n$, there exists an embedding $\psi$ of $T$ in $A\cup B$ with $\psi(r)=v$.
    Thus, this $A$ is as desired. 
\end{proof}

\subsection{Many leaves}
Here we will show the absorption for the many-leaves case. In this case, we will rely on the absorption of matching-makers.
\begin{lemma}\label{absorb:many leaves}
    Let $1/n\ll \varepsilon\ll \mu \ll \gamma$. Let $G$ be a $k$-graph with $\overline{\delta}_{\ell}(G)\geq \left(\overline{\delta}^{PM}_{k-1,\ell-1}+\gamma\right)$, and let $v\in V(G)$ be a vertex. Let $T$ be a hypertree on $\mu n$ vertices rooted at vertex $r$ with at least $2\mu^2 n$ leaf edges.

    Then there exists a vertex set $A\subseteq V(G)$ with $|A|=\left(\mu -\varepsilon\right)n$ such that for any $B\subseteq V(G)\setminus A$ with $|B|=\varepsilon n$, there exists an embedding $\psi$ of $T$ in $A\cup B$ with $\psi(r)=v$. 
\end{lemma}

\begin{proof}
    Let $\alpha>0$ be a constant such that $1/n\ll \alpha \ll \varepsilon$.
    Let $L$ be a set of exactly $\mu^2 n$ leaf edges of $T$ such that $r$ is not contained in any edge of $L$ as a leaf vertex.
    Let $T'=T\setminus L$. 
    Note that $T'$ is a hypertree with at most $(\mu -\mu^2)n$ vertices rooted at $r$. Since $\overline{\delta}_{d}\geq\left(\delta^{PM}_{k-1,d-1} +\gamma\right)$, we can greedily embed the vertices of $T'$ in~$G$ with $r$ mapped to $v$.
    Let the embedding function of $T'$ be $\psi'$ and $A'=\psi'(V(T'))$.
    Then we have $|A'|=(\mu -\mu^2)n$ and $\delta(G\setminus A')\geq \left(\delta^{PM}_{k-1,d-1} +\gamma/2\right)$.

    Next, we choose a suitable set $A''$ to embed the most of leaf edges. We randomly choose a vertex set $A''$ from $G\setminus A'$ by picking each vertex independently with probability $\frac{(k-1)\mu^2 n - \varepsilon n }{n-|A'|}$.
    Then by the Chernoff bound, $|A''| = (k-1)\mu^2 n - \varepsilon n \pm \alpha n$ with probability $1-o(1)$.
    In addition, by \Cref{lem:deg_concentration2}, with probability at least $1-o(1)$, for any $d$-tuple $S$ of $V(G)$, we have 
    \[
    d_{A''}(S)\geq \left(\delta^{PM}_{k-1,d-1} +\gamma/4\right)\binom{|A''|}{k-d}.
    \]
    By deleting or adding at most $\alpha n$ vertices to $A''$, we can ensure that $|A''|=(k-1)\mu^2 n - \varepsilon n$ and for any $d$-tuple $S$ of $V(G)$, 
    \[
    d_{A''}(S)\geq \left(\delta^{PM}_{k-1,d-1} +\gamma/5\right)\binom{|A''|}{k-d}.
    \]
    We claim that $A=A'\cup A''$ is a desired subset.
    We choose any $B\subseteq V(G)\setminus A$ with $|B|=\varepsilon n$ and let $L'=B\cup A''$. Then $|L'| = |A''| + |B| = (k-1)\mu^2 n$. As $\varepsilon\ll \gamma$, this implies that for any $d$-tuple $S$ of $V(G)$, we have
    \[
    d_{L'}(S)\geq \left(\overline{\delta}^{PM}_{k-1,d-1} +\gamma/5\right)\binom{|A''|}{k-d}\geq \left(\overline{\delta}^{PM}_{k-1,d-1} +\gamma/10\right)\binom{|L'|}{k-d}.
    \]
    Thus, by applying \Cref{lem:embed_stars} where $v_1, \ldots, v_m$ are the parent vertices of $L$ and $n_i$ is the number of leaves in $L$ adjacent to $v_i$, we can find a suitable embedding of $T$ in $A\cup B$ with $\psi(r)=v$.
    As $B$ is arbitrary, we can conclude that $A$ is a desired vertex set.
\end{proof}

By combining the two lemmas, we can prove the absorption lemma for hypertrees.
\begin{proof}[Proof of Lemma~\ref{absorb}]
    By \Cref{lem:leaf_or_bare_path}, $T$ contains either $\mu^2 n$ leaves or $2\mu^2 n$ edge-disjoint semi-bare paths of length $6$. The two cases follow from Lemma~\ref{absorb:many leaves} and Lemma~\ref{absorb:many bare paths}, respectively.
\end{proof}

\section{Proof of Theorem~\ref{main theorem}}\label{sec:main_theorem_proof}
We now prove our main theorem. 
\begin{proof}[Proof of Theorem~\ref{main theorem}]
    Let $\varepsilon, \mu>0$ be constants such that
    \[
    1/n \ll \eta \ll \varepsilon \ll \mu \ll \gamma, 1/k.
    \]
    By \Cref{lem:subtree_fixed_size}, $T$ can be decomposed into two subtrees $T_1, T_2$ such that $|V(T_1) \cap V(T_2)| =1$, $\mu n/2k \leq |V(T_1)| \leq 3k\mu n$ and $r \in V(T_2)$.
    Let $|V(T_1)| = \mu' n$. 
    Let $r'$ be the vertex in $V(T_1) \cap V(T_2)$ and $v' \in V(G) \setminus \{v\}$ be an arbitrary vertex.
    Then by \Cref{absorb}, there exists a vertex set $A \subseteq V(G)$ with $|A| = (\mu' - \varepsilon)n$ and $v\in A$ such that for any $B \subseteq V(G) \setminus A$ with $|B| = \varepsilon n$, there exists an embedding $\psi_1$ of $T_1$ in $A \cup B$ with $\psi_1(r') = v'$.
    As $\mu' \ll \gamma$, the minimum degree $\overline{\delta}_{\ell}(G \setminus A) \geq \left(\delta^{PM}_{k-1,\ell-1} + \gamma/2\right)$.
    By \Cref{lem:almost_spanning}, there exists an embedding $\psi_2$ of $T_2$ in $G \setminus A$ with $\psi_2(r) = v$ and $\psi_2(r') = v'$.
    Let $B = V(G) \setminus (A \cup \psi_2(V(T_2)))$.
    Then $|B| = \varepsilon n$ and thus there exists an embedding $\psi_1$ of $T_1$ in $A \cup B$ with $\psi_1(r') = v'$.
    By connecting the two embeddings $\psi_1$ and $\psi_2$, we can extend the embedding of $T$ to $G$ with $\psi(r) = v$.
\end{proof}

\section{Concluding remarks}\label{sec:concluding_remark}

In this paper, we determine the degree condition that ensures the existence of loose hypertrees.
In the study of hypertrees, many different types of hypertrees have been studied.
For instance, as we mentioned before, Pavez-Sign{'e}, Sanhueza-Matamala and~Stein~\cite{stein2020tree} discussed the tight hypertrees.
Recently, Rao, Sanhueza-Matamala, Sun, Wang and Zhou~\cite{rao2025degree} studied a class called expansion hypertrees, which is a special case of loose hypertrees.
Thus, we wonder whether the decomposition lemma (\Cref{lem:tree_split}) can be applied to study a different type of hypertrees.

We also note on the finding almost spanning hypertrees.
In~\cite{chen2025embedding}, the first author and Lo posed a problem about the degree condition that guarantees the existence of bounded degree almost spanning trees. We remark that in the proof of \Cref{lem:almost_spanning}, we only need the existence of a transversal almost-perfect matching instead of a transversal perfect matching.
Using this observation, we can show that if $\overline{\delta}_{k-1}> \frac{1}{k-1}$, an embedding for any almost spanning hypertrees exists. However, this bound may not be optimal for the bounded degree cases. Thus, we pose the following question.
\begin{question}
    Let $1/n \ll \varepsilon\ll \gamma \ll 1/\Delta, 1/k$. If $G$ is a
    $k$-graph on $n$ vertices with $\delta_{k-1}(G)\geq \left(1/k + \gamma\right)n$ and $T$ is a hypertree on $(1-\varepsilon)n$ vertices with $\Delta_1(T) \leq \Delta$, is there an embedding from $T$ to $G$?
\end{question}
\bibliographystyle{abbrv}
\bibliography{cite}

\appendix
\section{Proof of Lemma~\ref{lem:tree_split}}
We now prove Lemma~\ref{lem:tree_split}. The proof follows the argument in Im, Kim, Lee, and Methuku~\cite{Im2024} but we additionally need to care about the root vertex.

\begin{proof}[Proof of Lemma~\ref{lem:tree_split}]
    We construct a decreasing sequence of hypertrees $T = T'_0 \supseteq T'_1 \supseteq \cdots \supseteq T'_\ell$ and later let $T_i = T'_{\ell-i}$ for $i \in [\ell]$.
    Let $m = \lceil 10^3/\mu \rceil$. 
    We start with $T'_0 = T$. For any step $i \geq 0$, if $r_j$ is a leaf in $T'_i$ for some $j \in [2]$, then we let $e_{r_j}$ be the edge incident to $r_j$ and do not select it otherwise. 
    Assume that $T'_i$ has been defined for some $i \geq 0$ such that $r_1, r_2 \in V(T'_i)$.
    If $T'_i$ has a leaf matching of size at least $\mu n/(50 m (D+2))$, then we delete the edges in the leaf matching, excluding $e_{r_1}$ and $e_{r_2}$ if they are included, and let $T'_{i+1}$ be the resulting hypertree.
    We continue this process until we reach a hypertree $T'_k$ such that $T'_k$ does not have a leaf matching of size at least $\mu n/(50 m (D+2))$.
    Note that we delete at least $\mu n/(50 m (D+2))-2$ edges for each step so $k \leq 60 m D/\mu$. 
    For each parent vertex $v$ of $T'_k$, if $v$ is incident to at least $D+2$ leaf edges in $T'_k$, then we delete all the leaf edges incident to $v$ excluding $e_{r_j}$ if it is present.
    Let $S$ be the resulting hypertree.
    Then the leaf edges of $S$ are one of the following three types: (1) It contains a parent vertex of $T'_k$ as a leaf, (2) its parent vertex incident to at most $D+2$ leaf edges of $T'_k$, or (3) it is $e_{r_j}$ for some $j \in [2]$.
    In all three cases, the parent vertex of $T'_k$ corresponds to at most $D+2$ leaf edges in $S$.
    Therefore, the number of leaf edges in $S$ is at most $\mu n/50 m$.

    By \Cref{count bare path}, $S$ has a collection of edge-disjoint semi-bare paths $P_1, \ldots, P_s$ of length $m+1$ such that $e(S - P_1 - \cdots - P_s) \leq 3\mu n/25 + 2e(S)/(m+1) \leq \mu n/2$.
    As at most $\mu n/(50 m(D+2))$ vertices of $S$ have a different degree from $T'_k$, at least $\max\{0, s-\mu n/(50 m(D+2))\}$ of the semi-bare paths $P_i$ are still bare paths in $T'_k$.
    Let $Q_1, \ldots, Q_{s'}$ be the collection of such semi-bare paths of $T'_k$.
    If $r_j$ is an internal vertex of $Q_i$ for some $j \in [2]$, then we delete such a $Q_i$ from the list and reindex the remaining semi-bare paths as $Q_1, \ldots, Q_{s''}$.
    Each $Q_i$ contains a bare path of length $3$ in the middle. Let $T'_{k+1}$ be the forest obtained from $T'_k$ by removing these bare paths of length $3$.
    As each $Q_i$ is edge-disjoint, the deleted edges are vertex-disjoint from each other and we delete at most $n/m \leq \mu n$ paths in total.

    After constructing $T'_{k+1}$, we repeatedly delete one leaf edge from the remaining edges of each $Q_i$ except the edges at the end of each $Q_i$. 
    This gives $T'_{k+2}, \ldots, T'_{k+m-3}$ where each $T'_{k+i}$ is obtained by deleting a leaf matching from $T'_{k+i-1}$. In addition, $r_1, r_2 \in V(T'_{k+m-3})$.

    For each parent vertex $v$ of $T'_{k+m-3}$, if $v$ is incident to at least $D+2$ leaf edges in $T'_{k+m-3}$, then we delete all the leaf edges incident to $v$ excluding $e_{r_j}$ if it is included.
    Let $T'_{k+m-2}$ be the resulting hypertree and let $\ell = k+m-2$.
    Then $\ell \leq 60 m D/\mu + m - 2 \leq 10^5 D \mu^{-2}$.
    Thus, it suffices to show that $T'_\ell$ has at most $\mu n$ edges.
    If a vertex is a parent vertex of at least $D+2$ leaf edges in $T'_{k}$, then it still has at least $D+2$ leaf edges in $T'_{k+m-3}$. The way of obtaining $T'_{k+m-2}$ from $T'_{k+m-3}$ is the same as the way of obtaining $S$ from $T'_k$, so this implies that $E(T'_{k+m-2}) \subseteq E(S)$.
    As $T'_{k+m-3} = T'_k - Q_1 - \cdots - Q_{s''}$, we have that
    \begin{align*}
        e(T'_{k+m-2}) &\leq e(S-Q_1-\cdots-Q_{s''}) \leq e(S-P_1-\cdots-P_s) + m \cdot \left(\frac{\mu n}{50 m(D+2)}+2\right) \leq \mu n.
    \end{align*}
\end{proof}
\end{document}